\documentclass{article}
\usepackage[margin=1in]{geometry}
\usepackage{amsmath, amsfonts, amssymb, amsthm}
\usepackage[table]{xcolor}
\usepackage{mathtools}
\usepackage{array}
\usepackage{enumitem}
\usepackage{xcolor}
\usepackage{comment}
\usepackage[colorlinks=true, linkcolor=blue, citecolor=red, urlcolor=purple]{hyperref}

\usepackage{algorithm,algpseudocode}
\usepackage{authblk}
\usepackage{xspace}
\usepackage[capitalize,sort]{cleveref}
\usepackage{tikz}
\usepackage{thm-restate}
\usepackage{multirow}
\usepackage{multicol}
\usepackage[symbol]{footmisc}
\usepackage{physics}

\newtheorem{theorem}{Theorem}[section]
\newtheorem{conjecture}[theorem]{Conjecture}

\newtheorem{lemma}[theorem]{Lemma}
\newtheorem{observation}[theorem]{Observation}
\newtheorem{corollary}[theorem]{Corollary}
\newtheorem{definition}[theorem]{Definition}
\newtheorem{claim}[theorem]{Claim}
\newtheorem*{claim*}{Claim}
\newtheorem{example}[theorem]{Example}

\newtheorem{remark}[theorem]{Remark}
\newtheorem{fact}[theorem]{Fact}
\newenvironment{claimproof}[1]{\emph{Proof. }\space#1}{\hfill $\blacksquare$\medskip\par}

\newcommand{\I}{\mathcal{I}}
\newcommand{\M}{\mathcal{M}}
\newcommand{\B}{\mathcal{B}}
\newcommand{\C}{\mathcal{C}}
\newcommand{\A}{\mathcal{A}}
\newcommand{\floor}[1]{\lfloor #1 \rfloor}
\newcommand{\ceil}[1]{\lceil #1 \rceil}
\newcommand{\R}{\mathbb{R}}
\newcommand{\Rp}{\mathbb{R}_{\geq0}}

\DeclareMathOperator*{\Lex}{Lex}
\DeclareMathOperator*{\cl}{cl}

\title{Matroids are Equitable\footnote{An earlier version of this paper is to appear at SODA 2026 \cite{equitable-SODA}.}}
\date{}
\author[1,2]{Hannaneh Akrami}
\author[4]{Siyue Liu}
\author[3]{Roshan Raj} 
\author[1]{L\'aszl\'o A. V\'egh}
\affil[1]{Hertz Chair for Algorithms and Optimization, University of Bonn, Germany}
\affil[2]{Max Planck Institute for Informatics and Universität des Saarlandes}
\affil[3]{TIFR Mumbai}
\affil[4]{Tepper School of Business, Carnegie Mellon University, Pittsburgh, PA, USA}
\begin{document}
\maketitle
\thispagestyle{empty}
\begin{abstract}
We show that if the ground set of a matroid can be partitioned into $k\ge 2$ bases, then for any given subset $S$ of the ground set, there is a partition into $k$ bases such that the sizes of the intersections of the bases with $S$ may differ by at most one. This settles the matroid equitability conjecture by Fekete and Szab\'o (Electron.~J.~Comb.~2011) in the affirmative. We also investigate equitable splittings of two disjoint sets $S_1$ and $S_2$, and show that there is a partition into $k$ bases such that the sizes of the intersections with $S_1$ may differ by at most one and the sizes of the intersections with $S_2$ may differ by at most two; this is the best  one can hope for arbitrary matroids.

We also derive applications of this result into matroid constrained fair division problems. We show that there exists a matroid-constrained fair division that is envy-free up to one item if the valuations are identical and tri-valued additive. We also show that for bi-valued additive valuations, there exists a matroid-constrained allocation that provides everyone their maximin share. 
\end{abstract}

\pagenumbering{arabic}
\section{Introduction}
Let $\M=(E,\I)$ be a matroid defined on ground set $E$ with independent sets $\mathcal I$. We let $\B=\B(\M)$ denote the set of bases, i.e., maximal independent sets.\footnote{We include formal definitions of all concepts and notation in Section~\ref{sec:preliminaries}.} 
For $\ell\ge 1$, we say that $\M$ is \emph{$\ell$-equitable}, if the ground set $E$ is the disjoint union of two bases in $\M$, and the following property holds: for any $\ell$ disjoint sets $S_1,\ldots,S_\ell$, there exist two disjoint bases $B_1$ and $B_2$ such that $\lfloor |S_i|/2\rfloor \le |B_1\cap S_i| \le \lceil |S_i|/2\rceil$ for each $i\in [\ell]$.
We say that a class of matroids is $\ell$-equitable, if any matroid in the class, where the ground set is the disjoint union of two bases is $\ell$-equitable.
We say that $\M$ or a class of matroids is \emph{equitable} if it is 1-equitable.

This property was first investigated by Fekete and Szab\'o \cite{Fekete2011}. They showed that graphic matroids are equitable, and also showed 3-equitability of graphic matroids for the special case when $S_1$, $S_2$, and $S_3$ are disjoint stars. They also observed that weakly base orderable matroids are equitable, and that---using a result by Davies and McDiarmid \cite{Davies1976}---strongly base orderable matroids are $\ell$-equitable for any $\ell\ge 1$. A matroid is \emph{weakly base orderable} if for any bases $B_1$ and $B_2$, there is a bijection $f\, :\, B_1\to B_2$ such that $B_1-e+f(e)$ and $B_2-f(e)+e$ are bases for any $e\in B_1$,\footnote{We use  $X-x+y$ to denote $(X\setminus\{x\})\cup \{y\}$.}  and  \emph{strongly base orderable} if this bijection satisfies the stronger property that $(B_1\setminus X)\cup f(X)$ is a basis for any $X\subseteq B_1$.

The \emph{equitability conjecture} asks if all matroids are equitable. This was first posed in \cite{Fekete2011}, see also the website \cite{EgresOpen} in the Egres Open  collection. In recent years, progress has been made in multiple conjectures on matroid basis exchanges that  imply the equitability conjecture for certain classes. In particular, a result by B\'erczi, M\'atrav\"olgyi and Schwarcz \cite{Berczi2024} implies equitability of regular matroids. See Section~\ref{sec:related} for a discussion of these conjectures, as well as for other previously known special cases of the equitability conjecture.

Our main result resolves the equitability conjecture in the affirmative for all matroids. We state the result more generally, for partitions of the ground set into any number $k$ of bases. We note that equitability for two sets  easily implies this more general statement, see also the discussion page in \cite{EgresOpen}. Our result also includes the guarantee of polynomial-time computability.

\begin{theorem}{\label{thm:main-equitable}}
Let $\M=(E,\I)$ be a matroid, and assume that the ground set $E$ can be partitioned into the union of $k\ge 1$ disjoint bases. Let $S\subseteq E$ be an arbitrary set. Then, there exists a partition $E=B_1\cup B_2\cup\ldots\cup B_k$ into $k$ bases such that 
{\[
\left\lfloor \frac{|S|}{k}\right\rfloor \le |B_i\cap S|\le \left\lceil \frac{|S|}{k}\right\rceil\, \quad \forall i\in[k]\, .
\]}
Moreover, such bases can be found in polynomial time.
\end{theorem}

For $\ell\ge 2$, $\ell$-equitability fails in general. A simple example for $\ell=2$ is the graphic matroid of $K_4$, with $S_1$ and $S_2$ being two disjoint bases, see Example~\ref{exm:tight}. Nevertheless, we show near-equitability for $\ell=2$. Namely, equitability holds if at least one of $|S_1|$ and $|S_2|$ is odd, and is violated by one on either of the two sets if both have even cardinality. The example of $K_4$ demonstrates that this is the strongest statement one may hope for $\ell=2$. Our general statement is given  in Theorem \ref{thm:two-equitable} for $k$ bases; here, we present the corollary for two bases.

\begin{corollary}[of Theorem \ref{thm:two-equitable}]\label{cor:two-equitable}
    Let $\M=(E,\I)$ be a matroid, and assume that the ground set $E$ can be partitioned into the union of two disjoint bases. Let $S_1,S_2\subseteq E$ be two arbitrary disjoint sets. Then, there exists a partition $E=B_1\cup B_2$ into two bases such that 
    \begin{enumerate}[label = (\roman*)]
        \item If $|S_1|$ and $|S_2|$ are both odd, then $|B_1 \cap S_1| = \floor{|S_1|/2}$ and $|B_1 \cap S_2| = \ceil{|S_2|/2}$. %
        \item If $|S_1|$ is even and $|S_2|$ is odd, then $|B_1 \cap S_1| = |S_1|/2$ and $|B_1 \cap S_2| = \floor{|S_2|/2}$. 
        \item If $|S_1|$ and $|S_2|$ are both even, then $|B_1 \cap S_1| = |S_1|/2$ and $|B_1 \cap S_2| \in [|S_2|/2-1, |S_2|/2+1]$. 
    \end{enumerate}
\end{corollary}
Thus, in the third case, one of the even sets $S_1$ is partitioned equally, while the other, $S_2$ is partitioned such that $||S_2 \cap B_1| - |S_2 \cap B_2||\le 2$. Corollary~\ref{cor:two-equitable}
 proves a special case of a conjecture on almost fair representations in the intersection of matroids by Aharoni, Berger, Kotlar, and Ziv \cite{Aharoni2017}, as discussed in Section~\ref{sec:related}.

We next state the main technical tool used in proving the above results.
\begin{definition}
Let $B_1$ and $B_2$ be two disjoint bases in the matroid $\M=(E,\I)$. The set $X\subseteq E$ is an \emph{exchangeable set for $B_1$ and $B_2$}, if $X\subseteq B_1\cup B_2$, and  $B_1\Delta X$ and $B_2\Delta X$ are both bases.\footnote{$X\Delta Y$ denotes the symmetric difference of the sets $X$ and $Y$.} 
Further,
 for  $S\subseteq B_1\cup B_2$, and $t\in (B_1\cup B_2)\setminus S$ we say that $X\subseteq B_1\cup B_2$ is a \emph{$(t,S)$-exchangeable set for $B_1$ and $B_2$}, if $t\in X\subseteq S+t$ and $X$ is an exchangeable set for $B_1$ and $B_2$.
\end{definition}
Our next theorem shows the existence of an exchangeable set.
\begin{restatable}{theorem}{exchangethm}\label{thm:exchangable}
Let $\M=(E,\I)$ be a matroid, $S\subseteq E$, and let $B_1$ and $B_2$ be two disjoint bases. If $|B_1\cap S|<|B_2\cap S|$, then for some $t\in B_1\setminus S$ there exists a $(t,S)$-exchangeable set for $B_1$ and $B_2$. Moreover, such a set can be found in polynomial time.
\end{restatable}
 Theorem~\ref{thm:main-equitable} follows easily, as exchangeable sets can be used to `rebalance' the intersections of two bases $B_1$ and $B_2$ with a set $S$. That is, we repeatedly find $(t,S)$-exchangeable sets for bases $B_i$ and $B_j$ with $|B_i\cap S|<|B_j\cap S|-1$.

One may wonder if analogous results may exist for $\ell$-equitability. We formulate the following conjecture. 
\begin{conjecture}
     There exist a function $f\,:\,\mathbb{N}\to\mathbb{N}$ such that the following hold:
     For any matroid $\M=(E,\I)$ such that  the ground set $E$ can be partitioned into the union of $k$ disjoint bases, and any  disjoint subsets $S_1,S_2,\ldots,S_\ell$,  there exists a partition of $E$ into $k$ disjoint bases $B_1,\ldots,B_k$ such that $\left||B_i\cap S_j|-|S_j|/k\right\|\le f(\ell)$ for all $i\in [k]$ and $j\in[\ell]$. 
\end{conjecture} 

We note that the approach based on \Cref{thm:exchangable}
gives the above results for $\ell=1$ and $\ell=2$, but 
does not appear to yield a bound even for $\ell=3$. On the other hand, we are not aware of any example that shows $f(\ell)>2$. Thus, the constant function $f(\ell)=2$ could possibly satisfy the above conjecture.

\subsection{Applications in Fair Division}
We now state some applications of our results to discrete fair division settings.
A standard discrete fair division instance denoted by $\langle N,E,(v_i)_{i\in N}\rangle$ consists of a set $N=[n]$ of $n$ agents, a set $E$ of $m$ indivisible goods and 
a valuation function $v_i\,:\, 2^E\to \Rp$ for each agent $i\in N$. For ease of notation, we sometimes use $g$ instead of $\{g\}$ for good $g \in E$.
A feasible allocation is a subpartition $(A_1,\dots,A_n)$ of the set $E$, where $A_i$ denotes the bundle allocated to agent $i$, such that $\bigcup_{i \in N} A_i \subseteq E$ and $A_i \cap A_j = \emptyset$ for any $i \ne j$. While in a general case, any subpartition of $E$ into $n$ bundles is a feasible allocation, in many scenarios, certain restrictions may be imposed. For example, it might be undesirable to waste items and leave them unallocated, or items might be conflicting and a certain combination of the items cannot be in the same bundle.

A line of literature, starting with the paper by Biswas and Barman \cite{biswas2018fair}, examined matroid constrained allocations. Assume we are given a matroid $\M=(E,\I)$ on the set of goods. In a feasible allocation, it is required that {\em (i)}  all  goods must be allocated and {\em (ii)} all bundles must be independent sets of $\M$. The latter can be reduced to the case when each bundle is a basis \cite{Dror2023}. 
We denote a matroid constrained fair division instance as $\langle N,E,(v_i)_{i\in N},\M\rangle$. We formally define the allocations that we consider feasible in this work in Definition \ref{def:alloc}. 

\begin{definition}\label{def:alloc}
    Given a matroid constrained fair division instance $\langle N,E,(v_i)_{i\in N},\M\rangle$, ${A}=(A_1, \ldots, A_n)$ is a \emph{feasible allocation}, if
    \begin{enumerate}[label = (\roman*)]
        \item ${A}$ is a partition of $E$ into $n$ bundles, and
        \item for all $i \in N$, $A_i \in \B(\M)$ is a basis of $\M$.
    \end{enumerate}
\end{definition}

\begin{definition}[EF1]\label{def:ef1}
    An allocation ${A}$ is \emph{envy-free up to one good (EF1)}, if for all agents $i,j$ with $A_j \neq \emptyset$, there exists a good $g \in A_j$ such that $v_i(A_i) \geq v_i(A_j - g)$. 
\end{definition}

Given an allocation ${A}$, we say an agent $i$
\emph{envies} agent $j$ if  $v_i(A_i) < v_i(A_j)$, and $i$ {\emph{EF1-envies} agent $j$ if $v_i(A_i) < v_i(A_j-g)$ for all $g \in A_j$.
}
 
A valuation function $v\,:\,2^E\to \Rp$  is called \emph{additive} if $v(S)=\sum_{s\in g} v(g)$ for every $S\subseteq E$. Further,
an additive valuation function $v$ is \emph{binary} if for all $g \in E$, $v(g) \in \{0,1\}$. If all agents have identical binary valuations, Theorem \ref{thm:main-equitable} easily implies the existence of an EF1 allocation in the matroid constrained fair division problem \cite{Berczi2024}. It suffices to set $S$ to be the set of items of value $1$ for the agents. Then by Theorem \ref{thm:main-equitable}, one can (almost) equally partition $S$ into $n$ disjoint bases which results in an EF1 allocation if we allocate these bases to the agents arbitrarily.

\begin{corollary}[of Theorem \ref{thm:main-equitable}]
Given a matroid constrained fair division instance $\langle N,E,(v_i)_{i\in N},\M\rangle$, if $E$ can be partitioned into the union of $|N|$ disjoint bases, and the valuations $v_i$ are identical and binary additive, there exists a feasible EF1 allocation.
\end{corollary}

Extending the class of binary valuations, we
say that a valuation function $v\,:\,2^E\to\Rp$ is \emph{bi-valued} if there exists $a,b \in \Rp$ such that for all $g \in E$, $v(g) \in \{a,b\}$. Similarly, a valuation function $v$ is \emph{tri-valued} if there exists $a,b,c \in \Rp$ such that for all $g \in E$, $v(g) \in \{a,b,c\}$. Using Theorem~\ref{thm:exchangable}, we can also show the following.

\begin{restatable}{theorem}{trivalef}\label{thm:tri-val-ef1}
Given a matroid constrained fair division instance $\langle N,E,(v_i)_{i\in N},\M\rangle$, if $E$ can be partitioned into the union of $|N|$ disjoint bases, and the valuations $v_i$ are identical and tri-valued  additive, 
there exists a feasible EF1 allocation.
\end{restatable}

We also consider the \emph{maximin share fairness} notion.
The maximin share of an agent $i$ is the maximum value she can guarantee for herself if she could divide the goods into $n$ bundles and receive a minimum valued bundle. 
Let $\A$ denote the set of all feasible allocations. We adapt the notion of maximin share (MMS) to this setting and for a valuation function $v\,:\,2^E\to \Rp$ we define 
\[ \mu_{v}^n(E) \coloneqq \max_{A \in \A} \min_{j \in N} v(A_j)\, .\]

\begin{definition}
    An allocation $A$ is $\alpha$-MMS if and only if for all $i \in N$, $v_i(A_i) \geq \alpha \mu_{v_i}^n(E)$. For $\alpha=1$, we simply call the allocation MMS.
\end{definition}
From Theorem~\ref{thm:main-equitable}, we can infer the MMS values for binary additive valuations.
\begin{corollary}\label{cor:mms-val}
Let $\langle N,E,(v_i)_{i\in N},\M\rangle$ be  a matroid constrained fair division instance  such that the valuations $v_i$ are binary additive, and assume $E$ can be partitioned into the union of $n$ bases. Then, for all $i \in N$, $\mu_{v_i}^n(E) = \floor{v_i(E)/n}$. 
\end{corollary}
We show MMS allocations exist for the case of bi-valued additive valuations; note that these are not required to be identical.
\begin{restatable}{theorem}{mmsthm}\label{thm:MMS}
Given a matroid constrained fair division instance $\langle N,E,(v_i)_{i\in N},\M\rangle$, if $E$ can be partitioned into the union of $|N|$ disjoint bases, and all $v_i$'s are bi-valued additive valuations, then there exists a feasible MMS allocation.
\end{restatable}

\subsection{Related literature}\label{sec:related} 

\subsubsection{Matroid optimization}
Equitability is intimately related to multiple famous conjectures in matroid theory; we review three interrelated families of such conjectures.

\paragraph{Basis exchange sequences.}
Recall that by the strong bases exchange axiom of matroids, for any two bases $X,Y\in \B$, for any $e\in X\setminus Y$ there is an $f\in Y\setminus X$ such that $X-e+f,Y-f+e\in \B$.
White \cite{White1980} considered the question of exchanging basis sequences. Let ${\cal X}=(X_1,\ldots,X_k)$ and ${\cal Y}=(Y_1,\ldots,Y_k)$ be two collections of $k$ bases of a matroid $\M$. By a symmetric exchange from ${\cal X}$, we mean a symmetric exchange operation involving two bases in the collection. That is, for some $i,j\in [k]$, ${\cal X}'=(X_1,\ldots,X_{i-1},X'_i,X_{i+1},\ldots, X_{j-1},X'_j,X_{j+1},\ldots,X_k)$, where $X'_i=X_i-e+f$ and $X'_j=X_j-f+e$ for some $e\in X_i\setminus X_j$ and $f\in X_j\setminus X_i$. One may ask when the collection ${\cal Y}$ can be reached from ${\cal X}$ by a series of symmetric exchanges. A trivial necessary condition is that the two sequences are \emph{compatible}, i.e., each $e\in E$ appears in the same number of bases in both collections. White conjectured that this condition is also sufficient.
\begin{conjecture}[White, 1980]\label{conj:white}
The basis sequence ${\cal Y}$ can be reached from the basis sequence  ${\cal X}$ by a series of symmetric exchanges if and only if ${\cal X}$ and ${\cal Y}$ are compatible.
\end{conjecture}
The conjecture remains open even for $k=2$, and even if ${\cal X}=(B_1,B_2)$ and ${\cal Y}=(B_2,B_1)$, i.e., we need to swap two bases.
 For this setting, the following stronger conjecture was first formulated by Gabow \cite{Gabow1976}, see also \cite{Cordovil1993,Wiedemann}, and Conjecture 15.9.11 in Oxley's book \cite{Oxley2011}.
\begin{conjecture}[Gabow, 1976]\label{conj:gabow}
Given two bases $B_1=\{b_1,\ldots,b_r\}$ and $B_2=\{b'_1,\ldots,b'_r\}$ of a matroid $\M$, there exist permutations $\sigma$ and $\sigma'$ such that  in the sequence $(b_{\sigma(1)},\ldots,b_{\sigma(r)},b'_{\sigma'(1)},\ldots,b'_{\sigma'(r)})$, every set of $r$ cyclically consecutive elements forms a basis.
\end{conjecture}
Thus, this conjecture corresponds to a special case of White's conjecture, but also stipulates that in this case, the length of the exchange sequence is $r$, the rank of the matroid.

Observe that both these conjectures imply equitability. Indeed, assume that for any basis pair $(B_1,B_2)$, one  can reach $(B_2,B_1)$ through an (arbitrary long) series of symmetric exchanges. Without loss of generality, assume $|B_1\cap S|<|B_2\cap S|$. Since the size of the intersection $|B'_1\cap S|$ changes by at most one in each iteration, the first basis pair $(B'_1,B'_2)$ in the sequence with $|B'_1\cap S|\ge |S|/2$ will give an equitable partition. 
We note that Gabow's conjecture implies the existence of a matroid constrained EF1 allocation for identical additive valuations; we descibe the reduction in  \cref{lem:GabowImpliesEF1} in the Appendix.

In recent work B\'erczi, M\'atrav\"olgyi and Schwarcz \cite{Berczi2024} verified White's conjecture for $k=2$, as well as Gabow's conjecture for regular matroids. This in particular implies Theorem~\ref{thm:main-equitable} for regular matroids.  Quoting \cite[arXiv version]{Berczi2024}:
\begin{quote}\em
There are several further problems that aim at a better understanding of the structure of bases. The probably most appealing one is the Equitability Conjecture for matroids that provides a relaxation of both Conjecture 1 [White's conjecture] and Conjecture 2 [Gabow's conjecture].
\end{quote}
Further known cases of White's conjecture for $k=2$ include sparse paving matroids \cite{bonin2013basis}, split matroids \cite{berczi2024exchange}, strongly base orderable matroids \cite{lason2014toric}, and frame matroids satisfying a certain linearity condition \cite{mcguinness2020frame}. The papers \cite{bonin2013basis} and \cite{berczi2024exchange} also show Gabow's conjecture for the respective classes of sparse paving and split matroids. 
We refer the reader to \cite{Berczi2024} and Schwarcz's PhD thesis \cite{SchwarczPhD} for further background and references on these conjectures.

We note that while Conjectures~\ref{conj:white} and \ref{conj:gabow} immediately yield 1-equitability, they do not appear to imply any variant of 2-equitability. In particular, we do not see a way how Corollary~\ref{cor:two-equitable}, or our main technical tool, Theorem~\ref{thm:exchangable} on $(t,S)$-exchangeable sets would follow.

\paragraph{Partitioning into common independent sets.}
The following natural question in matroid optimization unifies several known results, as well as includes multiple longstanding conjectures: 

\begin{quote}
Given two matroids $\M_1=(E,\I_1)$ and $\M_2=(E,\I_2)$ on a common ground set $E$ and the same rank $r$, such that the ground set can be partitioned into $k$ bases both in $\M_1$ and in $\M_2$, what is the minimum number of common independent sets covering $E$?\end{quote}

The most interesting case is when the answer is $k$, i.e., there exists a partition of the ground set into $k$ common bases. Classical examples include edge colourings of bipartite graphs and Edmonds's disjoint arborescences theorem.  B\'erczi and Schwarcz \cite{berczi2021complexity} proved the exponential oracle complexity of deciding the existence of $k$ disjoint common bases; we refer the reader to this paper for background and known results. 
$\ell$-Equitability can be phrased as a special case of this question where $k=2$, the first matroid has rank $r$, and the second matroid is a partition matroid.\footnote{Let us define $T_i=S_i$ if $|S_i|$ is even and $T_i\subseteq S_i$ with $|T_i|=|S_i|-1$ if $|S_i|$ is odd; let $R=E\setminus(\cup_i T_i)$. We consider the partition matroid for the partition $(T_1,\ldots,T_k,R)$, with $r(T_i)=|T_i|/2$ and $r(R)=r-\sum_i |T_i|/2$.  It is easy to see that $\ell$-equitability is equivalent to partitioning the ground set into two common bases of these matroids.}

Davies and McDiarmid \cite{Davies1976} proved that if both matroids are strongly base orderable, and the ground set can be partitioned into $k$ bases in each of the matroids, then it can also be partitioned into $k$ common bases. Since partition matroids are strongly base orderable, this implies $\ell$-equitability for any strongly base orderable matroid.

An intriguing conjecture of the above form is Rota's basis conjecture: given a matroid of rank $r$ whose ground set partitions into $r$ disjoint bases $B_1,B_2,\ldots,B_r$, there always exists a transversal set of disjoint bases, i.e., $B'_1,B'_2,\ldots,B'_r$ such that $|B_i\cap B'_j|=1$ for all $i,j\in[r]$. Note that this can be seen as a special instance of $r$-equitability.

Regarding the general question, Aharoni and Berger \cite{Aharoni2006} conjectured that if the ground set can be partitioned into $k$ disjoint bases in both matroids, then it can be covered by at most $k+1$ common independent sets. They showed the weaker bound that there always exists a partition into $2k$ common independent sets. This in particular shows that in Rota's conjecture, the ground set can be covered by $\le 2r$ independent sets that intersect each $B_i$ at most once. On the other hand, Buci{\'c},  Kwan, Pokrovskiy, and Sudakov \cite{bucic2020} showed the existence of $(\frac{1}2-o(1))n$ disjoint transversal bases, and 
Pokrovskiy \cite{pokrovskiy2020rota} gave an asymptotic version showing the existence of $(1-o(1))n$ disjoint transversal independent sets of size  $(1-o(1))n$. Recently, Montgomery and Sauermann \cite{montgomery2025} improved these results by proving the existence of $(1-o(1))n$ disjoint transversal bases.

\paragraph{Fair representations of matroid intersections.}
A related notion of \emph{fair representations} was introduced by Aharoni, Berger, Kotlar, and Ziv \cite{Aharoni2017}. Let us say that for sets $S,X\subseteq E$, and $\alpha\in[0,1]$, $X$ represents $S$ 
\emph{$\alpha$-fairly}, if $|S\cap X|\ge \lfloor \alpha|S|\rfloor$, and 
\emph{almost $\alpha$-fairly}, if $|S\cap X|\ge \lfloor \alpha|S|\rfloor-1$. $X$ represents a partition of $E$ (almost) $\alpha$-fairly, if it represents all of the sets  (almost) $\alpha$-fairly.
They formulated the following general conjecture. 
\begin{conjecture}[Aharoni, Berger, Kotlar, and Ziv, 2017]\label{conj:almost-fair}
Let $\M_1$ and $\M_2$ be two matroids on the same ground set $E$, and assume $E$ can be partitioned into $k$ independent sets in both matroids. For any partition of $E$, there exists a common independent set $X$ in the two matroids that represents this partition almost $1/k$-fairly.
\end{conjecture}
Assume $\ell$-equitability holds for a matroid $\M$ such that the ground set is the disjoint union of two bases of $\M$. Then, it follows that for every partition of the ground set with $\ell+1$ parts, there exists a common independent set of $\M$ and its dual matroid $\M^\star$ that represents the partition almost $1/2$-fairly. Indeed, let $E=S_0\cup S_1\cup\ldots\cup S_\ell$ be a partition into $\ell+1$ sets. By $\ell$-equitability, there exists two disjoint bases $B_1$ and $B_2$ with $\lfloor |S_i|/2\rfloor \le |B_1\cap S_i|,|B_2\cap S_i|$. Assume without loss of generality that $|B_1\cap S_0|\ge |B_2\cap S_0|$. Then, $B_1$ gives a fair representation of the partition, and is a basis both in $\M$ and in $\M^\star$.

In \cite{Aharoni2017}, it was noted the graphic matroid of $K_4$ shows that the stronger notion of fair representation is not possible. They also proved that if the common ground set $E$ can be partitioned into two independent sets in two matroids $\M_1$ and $\M_2$, then  any set $S\subseteq E$ has an almost $\left(\frac{1}{2}-\frac{1}{|E|}\right)$-fair representation.

In this context, Theorem~\ref{thm:main-equitable} and Corollary~\ref{cor:two-equitable} verify Conjecture~\ref{conj:almost-fair} for the case when the two matroids are each others duals, $\M$ and $\M^\star$,  $E$ is the disjoint union of two bases of $\M$, and the partition has at most three parts.

\paragraph{Exact Matroid Intersection.} Given a graph with the edges colored either red or blue and an integer $k$, the \emph{exact perfect matching problem} asks to find a perfect matching with exactly $k$ red edges or to determine that such a matching does not exist. Mulmuley, Vazirani, and Vazirani \cite{MVV87} gave an efficient randomized algorithm for this problem via the Isolation Lemma. However, it is a long standing open question to design a deterministic polynomial time algorithm even for the bipartite case. 

The \emph{Exact Matroid Intersection} problem generalizes the exact perfect matching problem for bipartite graphs. Here, we are given an integer $k$ and two matroids of same rank on a ground set whose elements are colored either red and blue. The goal is to find a common basis with exactly $k$ red elements. 
Exact matroid intersection turns out to be hard in general: H\"orsch et al.~\cite{Hoersch_Imolay_Mizutani_Oki_Schwarcz_ICALP2024}~proved that an easier task of finding a common basis with an odd (even) number of red elements requires exponentially many independence queries. 
An efficient randomized algorithm
is known when the input matroids are linear \cite{NSV94}. Deterministic polynomial-time algorithms are known for Pfaffian pairs \cite{webb2004counting}, for arborescences (that are intersections of graphic and partition matroids) \cite{DBLP:journals/tcs/ItoIKNOW23}, and for union of arborescences as well \cite{kobayashi2025reconfiguration}.

In this context, we get a deterministic polynomial time algorithm for exact matroid intersection in the special setting when the input matroids are dual of each other. We can first find a common basis with minimum number of red elements using weighted matroid intersection. Then, using \cref{thm:exchangable} iteratively, we can get a common basis with exactly $k$ red elements if such a common basisexists.

\subsubsection{Fair division}
Fair division of indivisible items has gained significant attention in the last decade. In the unconstrained setting, EF1 exists and can be 
computed in polynomial time when agents have monotone valuations whether the set of items consists of goods \cite{Lipton2004} or chores \cite{Bhaskar2021}.\footnote{
 In the general case allowing both goods and chores, when $v_i(g)$ could be positive, negative or zero, an allocation ${A}$ is \emph{EF1}, if for all agents $i,j$, there exists an item $g \in A_i \cup A_j$ such that $v_i(A_i - g) \geq v_i(A_j - g)$, unless $A_i=A_j=\emptyset$. Note that for the goods case this coincides with \Cref{def:ef1}.}
On the other hand, MMS allocations are not guaranteed to exist (even) when agents have additive valuations \cite{feige2021tight,Kurokawa2018}, and hence the focus has moved to approximations of MMS; for instance, guaranteeing $\alpha$ fraction of the MMS values of all the agents.  The state-of-the-art for additive valuations is $\alpha = 3/4 + 3/3836$ for goods \cite{DBLP:conf/soda/AkramiG24} and $\alpha = 13/11$ for chores \cite{Huang2023}. We refer the reader to the survey \cite{AMANATIDIS2023103965} for a more detailed discussion on the literature. For the rest of the overview, we limit our focus to \textit{constrained fair division}.

In constrained fair division, instances are furnished by a set of constraints on which combination of items can be allocated to the same agent. Li and Vetta \cite{Li2021} initiated the study of approximating MMS allocations under hereditary set systems. In this set system, any subset of an independent set is also independent and the agents can only receive independent sets. Hummel \cite{Hummel2025} improved the approximation factor in this setting to $2/3$. Hereditary set system valuations are a subclass of XOS valuations. A series of work on approximating MMS for XOS valuations \cite{GHODSI2022103633,SEDDIGHIN2024104049} led to the state-of-the-art of $3/13$ \cite{Akrami2023}. MMS has also been studied for other classes of valuation functions such as submodular \cite{Barman2020,GHODSI2022103633,uziahu2023} and subadditive \cite{feige2025,SEDDIGHIN2024104049,seddighin2025}. Gourv\`es and Monnot \cite{GOURVES201950} studied MMS under matroid constraints and proved $1/2$-MMS allocations can be computed in polynomial time.

Barman and Biswas \cite{biswas2018fair} studied fair allocation under partition matroid constraint and established efficient algorithms to compute EF1 and approximately MMS allocations. They also proved that when agents have identical additive valuations, under laminar matroid constraint, EF1 allocations can be computed efficiently.\footnote{This paper, as well as the subsequent  \cite{Biswas_Barman_2019}, claimed existence and efficient computability of EF1 allocations under general matroid constraints. However, the key exchange lemma in these conference papers was flawed, see the arxiv version of \cite{biswas2018fair}.} It is also known that allocations which maximize the geometric mean of agents' utilities (i.e. Nash welfare) are $1/2$-EF1 and Pareto optimal under arbitrary matroid constraints \cite{Cookson2025,Wang2024}.

However, very little is known regarding EF1 without approximation when agents have non-identical valuations, or non-identical matroid constraints, or when the underlying matroid is an arbitrary matroid. Dror, Feldman, and Segel-Halevi \cite{Dror2023} presented polynomial algorithms to compute EF1 allocations for arbitrary number of agents with heterogeneous binary valuations and partition matroids with heterogeneous capacities. They also proved stronger results for the case of two or three agents. In particular, they considered identical base orderable matroid constraints and showed the existence of EF1 allocation when either there are three agents with binary additive valuations or there are two agents with general additive valuations. Using the cut-and-choose method, \cref{thm:main-equitable} implies the existence of an EF1 allocation when one of the agents has a binary additive valuation.

To the best of our knowledge, this is the first work establishing fairness guarantees for general matroids; see \Cref{thm:tri-val-ef1,thm:MMS}. The existence of EF1 allocations for general matroids under additive valuations, even when agents have identical valuations, remains an open question.
 
We refer the reader to the surveys of \cite{Suksompong2021} and \cite{Biswas2023} for more discussion on constrained fair allocation.

\subsection{Overview of techniques}
The proof of Theorem~\ref{thm:exchangable} is based on a structural analysis of the \emph{exchange graph} $D(B_1,B_2)$ between $B_1$ and $B_2$ (see Section~\ref{sec:cycles_and_exchangeable_sets}). There is an edge $(x,y)$ from $x\in B_i$ to $y\in B_{3-i}$ if $B_i-x+y$ is a basis. Assume $|B_1 \cap S| < |B_2 \cap S|$. 
The goal is to find an exchangeable set $X$ such that $|(B_1 \Delta X) \cap S| = |B_1 \cap S|+1$. 
Pairs where $(x,y)$ and $(y,x)$ are both present represent symmetric exchanges. If there is a symmetric  edge $(t,y)$ with $t\in B_1\setminus S$ and $y\in B_2\cap S$, then $X=\{t,y\}$ is a $(t,S)$-exchangeable set and we are done. This is in fact always the case when the matroid is 
base orderable. Such matroids have a perfect matching of symmetric exchange edges; by the assumption $|B_1\cap S|<|B_2\cap S|$, a perfect matching must include a symmetric exchange edge $(t,y)$ as above.

The above argument however fails  for general matroids: while for each $x\in B_1\cup B_2$, there is a symmetric exchangeable element, it could be the case that the elements $x\in S\cap B_2$ can only be exchanged to elements in $S\cap B_1$. Instead, we show the existence of a directed cycle $\C$ in the exchange graph such that $V(\C)\setminus S=\{t\}$ for some $t\in B_1$. Using standard basis exchange arguments, the existence of such a $\C$ implies the existence of a $(t,S)$-exchangeable set (Lemma~\ref{lem:CycleImpliesExchangibility}). 
This proof also naturally gives rise to a polynomial-time algorithm: the main task is to construct the exchange graph, and find a chordless cycle $\C$ containing some $t\in B_1\setminus S$ with $V(\C)\subseteq S+t$. 

To prove the existence of such a directed cycle $\C$, we consider the restriction of the exchange graph to $S$.  Pick a strongly connected component $K\subseteq S$ with $|K\cap B_1|<|K\cap B_2|$. We crucially argue that for some $t\in B_1\setminus K$, there exist two edges $(x,t),(t,y)$ for some $x,y\in B_2\cap K$. This in particular implies $t\notin S$ since otherwise $t$ is also in the same strongly connected component $K$. Furthermore, there is a directed path from $y$ to $x$ inside $K$, which together with the two edges $(x,t),(t,y)$ forms a directed cycle $\C$ satisfying $V(\C)\setminus S=\{t\}$.

As already noted, Theorem~\ref{thm:main-equitable} on equitability follows easily from Theorem~\ref{thm:exchangable}. To prove Corollary~\ref{cor:two-equitable} on 2-equitability with respect to two sets $S_1$ and $S_2$, and the more general form Theorem~\ref{thm:two-equitable} in Section~\ref{sec:two-equitable}, we consider a partition into disjoint bases that minimizes certain potentials capturing balancedness. One can argue that if the near-2-equitability statements do not hold, then one can apply Theorem~\ref{thm:exchangable} for some $S\in \{S_1,S_2,E\setminus (S_1\cup S_2)\}$ to find a better partition.

Section~\ref{sec:fair-division} derives the applications to fair division problems. In particular, Theorem~\ref{thm:tri-val-ef1} shows the existence of matroid constrained EF1 allocations for identical tri-valued additive valuations. First, we observe that without loss of generality we can assume that the potential utility values are $0<a<b$ (Lemma~\ref{lem:shift-vals}). We use the near-2-equitability partition with respect to the set of items of values $b$ and $a$. This will already imply EF1 if $b\ge 2a$. However, if $b<2a$, additional exchanges are needed to transform the allocation to EF1.

Finally, the proof of Theorem~\ref{thm:MMS} uses the \emph{lone divider} technique \cite{AIGNERHOREV2022164}. In this approach, in order to guarantee a threshold of $t_i$ for each agent $i$, we inductively assign bundles to agents and eliminate them. In each step, a remaining agent $i$ (the lone divider) partitions the set of the remaining items into $n'$ bundles of value at least $t_i$ to herself, where $n'$ is the number of the remaining agents. Then we allocate a subset of these bundles to a subset of agents and eliminate them such that {\em(i)} all the eliminated agents are satisfied with their bundles, and {\em(ii)} all the remaining agents do not desire any of the allocated bundles. When using this approach, the main part of the proof is to show that a lone divider exists in each iteration \cite{Akrami_Rathi_2025-2,Akrami_Rathi_2025,amanatidis2021maximum,Hummel2025}. The rest follows by carefully applying Hall's theorem.
 
 \section{Preliminaries}\label{sec:preliminaries}
For a set $X$ and an element $e$, we use $X-e$ and $X+e$ as shorthand notations for $X \setminus \{e\}$ and $X \cup \{e\}$ respectively. We let $X\Delta Y$ denote the symmetric difference of the sets $X$ and $Y$, and $X\sqcup Y$ their disjoint union. For a directed graph $D=(V,A)$ and $S\subseteq V$, we let $A[S]=\{(i,j)\in A\mid i,j\in S\}$ denote the restriction of the edge set to $S$.

 A \emph{matroid} $\M$ is a set family $(E,\I)$ defined on \emph{ground set} $E$ such that $\I$ is a family of subsets of $E$ satisfying the following properties: 
	 \begin{enumerate}[label=(\roman*)]
        \item $\emptyset \in \I$.
		\item \emph{Downward-closed property:} For every $I\in \I$, all its subsets are also in $\I$.
		\item \emph{Augmentation property:} For every $I, J\in \I$ with $|J|>|I|$, there exists an element $a\in J\setminus I$ such that $I\cup\{a\}\in \I$.
	\end{enumerate}
	The sets in $\I$ are called the \emph{independent} sets of $\M$; all other sets are called \emph{dependent}. An independent set is called a \emph{base} if it is an inclusion-wise maximal set in $\I$; according to {\em (iii)}, these all have the same cardinality which we call the \emph{rank} of the matroid. Similarly, the rank of a set $S \subseteq E$ is the size of the maximum independent subset of $S$ and is denoted by $r(S)$.

The \emph{closure} (or \emph{span}) of a set $S \subseteq E$, is defined as $\cl(S) = \{e \in E \mid r(S+e)=r(S)\}$. We summarize some simple properties of the closure.
\begin{fact}\label{fact:closure}
In a matroid $\M=(E,\I)$,
\begin{enumerate}[label=(\roman*)]
    \item $S\subseteq \cl(S)$ for every $S\subseteq E$;
    \item $\cl(S)\subseteq \cl(T)$ for every $S\subseteq T\subseteq E$;
    \item $\cl(\cl(S))=\cl(S)$ for every $S\subseteq E$.
\end{enumerate}
\end{fact}
For every basis
$B\in \B$ and $e\in E\setminus B$, $B+e$ contains a unique circuit. 
 This is called the \textit{fundamental circuit} of $e$ with respect to $B$, and is denoted by $C(B,e)$. We state this in the following lemma, with an additional characterization:
 \begin{lemma}
     Let $\M=(E,\I)$ be a matroid and $B$ a basis, $e\in E\setminus B$. Then, $B+e$ contains a unique circuit denoted as $C(B,e)$, and
     \[
     C(B,e)=\{f\in B+e \mid B-f+e \in \B\}\, .
     \]
 \end{lemma}

For an independent set $I\in\I$, let $D(I)$ denote a directed bipartite graph with vertex set $E=I\sqcup (E\setminus I)$ and edge set $$\{(x,y)\mid x\in I, y\in E\setminus I, I-x+y \in \mathcal I\}.$$ 
Thus, for a basis $B$ and $e\in E\setminus B$, $C(B,e)$ is the set formed by $e$ and its in-neighbours in $D(B)$.
The following claim follows immediately.
\begin{observation}\label{obs:span}
    For a basis $B$ and $e\in E\setminus B$, let $Z$ be the set of in-neighbours of $e$ in $D(B)$. Then $e \in \cl(Z)$. 
\end{observation}

\begin{theorem}[{\cite[Corollary 39.12a, Theorem 39.13]{Sch03B}}]\label{thm:matchingExchange}
Let $I$ be an independent set in the matroid $\M=(E,\I)$, and $J$ another set with $|I|=|J|$. Then, \begin{itemize}
        \item If $J$ is an independent set, then $D(I)[I\Delta J]$, i.e., the subgraph of $D(I)$ with vertex set $I\Delta J$,  contains a perfect matching.
        \item If $D(I)[I\Delta J]$ has a unique perfect matching, then $J$ is also an independent set.
    \end{itemize}
\end{theorem}

\section{Proof of the Equitability Theorem}\label{sec:main}
In this section, we prove the main exchangeability result, Theorem~\ref{thm:exchangable}, and derive Theorem~\ref{thm:main-equitable}. First, in Section~\ref{sec:cycles_and_exchangeable_sets}, we start by showing how cycles in the exchange graphs give rise to exchangeable sets.
\subsection{Cycles and Exchangeable Sets}\label{sec:cycles_and_exchangeable_sets}
For two disjoint bases $B_1$ and $B_2$, let $D(B_1, B_2)$ be the exchange graph between $B_1$ and $B_2$, where an edge $(x,y)$ from $x\in B_i$ to $y\in B_{3-i}$ exists if $B_i-x+y$ is a basis. In other words, it is the directed graph on vertex set $B_1\cup B_2$ and the edges of $D(B_1)$ and $D(B_2)$ restricted to this set.
By a \emph{cycle}, we mean a directed cycle in this graph. For a cycle $\C$, we let $V(\C)$ denote the vertex set of $\C$.
By a \emph{chordless cycle}, we mean a cycle such that $V(\C)$ does not include any other edges than $\C$.
\begin{lemma}\label{lem:minimalCycleExchange} Let $B_1$ and $B_2$ be two disjoint bases and $\C$ be a chordless cycle in graph $D(B_1,B_2)$. Then, $V(\C)$ is an exchangeable set for $B_1$ and $B_2$.
\end{lemma}
\begin{proof}
Let $U=V(\C)$.
    Since  $\C$ is a chordless cyle  in $D(B_1,B_2)$, both $D(B_1)[U]$ and $D(B_2)[U]$ have unique perfect matchings. Hence, from \cref{thm:matchingExchange}, $B_1\Delta U$ and $B_2\Delta U$ are both bases.
\end{proof}
As the next lemma shows, cycles with chords also guarantee the existence of exchangeable sets. 
\begin{lemma} \label{lem:CycleImpliesExchangibility}
    Let $B_1$ and $B_2$ be two disjoint bases and $\C$ be a cycle in graph $D(B_1,B_2)$; let   $t\in V(\C)$. Then, there exists an exchangeable set $U\subseteq V(\C)$ for $B_1$ and $B_2$ with $t\in U$. Moreover, if $\C$ is such that no cycle $\C'$ exists with  $t\in V(\C')\subsetneq V(\C)$, then $V(\C)$ itself is exchangeable.
\end{lemma}
\begin{proof}
Let $A$ denote the edge set of $D(B_1,B_2)$. Without loss of generality, assume $t\in B_1$. To prove both parts of the statement,
    let $\C'$ be a cycle with minimum number of vertices such that $t\in V(\C')\subseteq V(\C)$. Note that $\C'$ might have a chord. We show that $U=V(\C')$ is an exchangeable set.
    Let $\C'=(v_0,v_1,v_2,\dots,v_{2k-2},v_{2k-1})$ for $v_0=t$. By the minimal choice of $\C'$, there is no edge $(v_i,v_j)\in A$ with $j>i+1$ or $i<2k-1$ and $j=0$, as such an edge would lead to a shorter cycle. 
    
    The edges $(v_0,v_1),(v_2,v_3),\dots$ $, (v_{2k-2},v_{2k-1})$ form a perfect matching $N_1$ from $B_1\cap U$ to $B_2\cap U$, and $(v_1,v_2)$, $(v_3,v_4),\ldots,(v_{2k-3},v_{2k-2}),(v_{2k-1},v_0)$ form a perfect matching $N_2$ from $B_2\cap U$ to $B_1\cap U$. We show that both are unique perfect matchings from $B_1\cap U$ to $B_2\cap U$ and from $B_2\cap U$ to $B_1\cap U$, respectively. 
    Then by \cref{thm:matchingExchange}, $U$ is an exchangeable set.

    Indeed, let $N'_1\neq N_1$ be another perfect matching from $B_1\cap U$ to $B_2\cap U$. Let $0\le i\le k-1$ be the smallest index such that $v_{2i}$ is not matched to  $v_{2i+1}$ in $N_1'$.
    Then, we must have $(v_{2i},v_j)\in A$ for some $j>2i+1$, but no such edge may exist as shown above. Thus, no such $N'_1$ may exist.

     Next, let $N'_2\neq N_2$ be another perfect matching from $B_2\cap U$ to $B_1\cap U$. Let $0\le i\le k-2$ be the smallest index such that $v_{2i+1}$ is not matched to  $v_{2i+2}$ in $N_2'$. (If no such $i$ exists, then $N'_2=N_2$ since $v_{2k-1}$ must also be matched to the remaining node $v_0$.) Again, we must have $(v_{2i+1},v_j)\in A$ for some $j>2i+2$ or for $j=0$, both giving a contradiction. Hence, no such $N'_2$ may exist either, completing the proof.
    \end{proof}

\subsection{Finding an exchangeable set}
We now prove \Cref{thm:exchangable}, restated for convenience.
\exchangethm*
\begin{figure}
    \centering
    \includegraphics[width=0.32\linewidth]{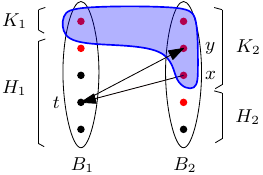}
    \caption{The elements in $S$ are colored red. The blue set illustrates $K$.}
    \label{fig:proof}
\end{figure}
\begin{proof}
We consider the directed graph $D(B_1,B_2)$; let $A$ denote its edge set.
Since $|B_1\cap S|<|B_2\cap S|$, by averaging there exists a strongly connected component $K\subseteq S$ in $A[S]$ such that $|K \cap B_1|<|K\cap B_2|$. 
Let $K_i=B_i\cap K$ and $H_i = B_i \setminus K_i$. Our goal is to  prove that there exists $t \in H_1$ such that edges $(x,t)$ and $(t,y)$ exist for some $x,y \in K_2$; see Figure \ref{fig:proof}. Since $|H_1| > |H_2|$ and both $H_1$, $H_2$ are independent, $H_1 \setminus \cl(H_2) \neq \emptyset$.

\begin{claim}\label{claim:1}
    For every $t \in H_1 \setminus \cl(H_2)$, there exists an edge $(x,t)$ for some $x \in K_2$.
\end{claim}
\begin{claimproof}
    Towards contradiction, assume for some $t \in H_1 \setminus \cl(H_2)$, all its in-neighbors are in $B_2 \setminus K_2 = H_2$. By \Cref{obs:span}, $t \in \cl(H_2)$,  a contradiction.
\end{claimproof}

\begin{claim}\label{claim:2}
    There exists a $(t,y)$ edge with $t \in H_1 \setminus \cl(H_2)$ and $y \in K_2$.
\end{claim}
\begin{claimproof}
    Towards contradiction, assume all in-neighbors of $K_2$ are in $B_1 \setminus (H_1 \setminus \cl(H_2))=K_1 \cup (\cl(H_2) \cap B_1)$. By \Cref{obs:span}, $K_2 \subseteq \cl(K_1 \cup (\cl(H_2) \cap B_1))$. 
    Therefore, 
\begin{align*}
    B_2 &= K_2 \cup H_2 \\
    &\subseteq \cl(K_1 \cup (\cl(H_2) \cap B_1)) \cup \cl(H_2) \\
    &\subseteq \cl(K_1 \cup (\cl(H_2) \cap B_1) \cup H_2) \tag{$\cl(S) \cup \cl(T) \subseteq \cl(S \cup T)$ for all $S,T \subseteq E$} \\
    &\subseteq \cl(K_1 \cup \cl(H_2))\\
    &=\cl(K_1\cup H_2). 
\end{align*}
The last equality holds because $\cl(S\cup T)\subseteq \cl(S\cup \cl(T))\subseteq \cl(\cl(S\cup T))=\cl(S\cup T)$ for all $S,T\subseteq E$ (see Fact~\ref{fact:closure})
This implies,
$$|B_2| = r(B_2) \leq r(K_1 \cup H_2) \leq  r(K_1) + r(H_2) = |K_1| + |H_2| < |K_2|+|H_2|=|B_2|,$$ 
which is a contradiction.
\end{claimproof}

Claims \ref{claim:1} and \ref{claim:2} imply that there exits $t \in B_1 \setminus K$ and $x,y \in K \cap B_2$, such that the edges $(x,t)$ and $(t,y)$ exist. Since $K$ is a strongly connected component of $A[S]$, if $t \in S$, we must have $t \in K$. Therefore, $t \in B_1 \setminus S$. By strong connectivity of $K$, there exists a path from $y$ to $x$ in $A[K]$. Hence there exists a cycle $\C$ with $t \in V(\C) \subseteq K+t \subseteq S+t$. By \Cref{lem:CycleImpliesExchangibility}, a $(t,S)$-exchangeable set exists.

Regarding polynomial time computability, \Cref{lem:CycleImpliesExchangibility} shows that if $\C$ is a cycle with $|V(\C)\setminus S|=1$, $V(\C)\setminus S\subseteq B_1$, and $V(\C)$ is minimal subject to this property, then $V(\C)$ is exchangeable. Moreover, the proof above guarantees the existence of such a cycle. We can consider any $t\in B_1\setminus S$, and solve a simple reachability problem to see if there is a cycle in $A[S+t]$ containing $t$, and find a shortest such cycle if one exists. For some $t$, we are guaranteed to find a $(t,S)$-exchangeable set this way.
\end{proof} 

\begin{proof}[Proof of \Cref{thm:main-equitable}]
We start from any partition of $E$ into $k$ disjoint bases $B_1,B_2,\ldots,B_k$. We execute the following algorithm. In each iteration, let $B_i$ be a basis with $|B_i\cap S|$ minimal, and $B_j$ a basis with $|B_j\cap S|$ maximal. If $|B_i\cap S|\le |B_j\cap S|-2$, then we use \Cref{thm:exchangable} for $B_i$ and $B_j$ to find $B'_i$ and $B'_j$ with $B'_i\cup B'_j=B_i\cup B_j$ and $|B'_i\cap S|=|B_i\cap S|+1$. This process clearly terminates in $O(|S|)$ steps with a partition such that $\max_i |B_i\cap S|-\min_i |B_i\cap S|\le 1$.
\end{proof}

\section{Equitability for Two Sets}\label{sec:two-equitable}
In this section, we consider $2$-equitability, i.e., partitioning the ground set of a matroid $\M=(E,\I)$  into $k$ disjoint bases such that each of them intersect two given disjoint sets $S_1$ and $S_2$ as equally as possible. Theorem~\ref{thm:two-equitable} below is our general statement for $k$ bases; for $k=2$, we obtain  Corollary~\ref{cor:two-equitable} stated  in the Introduction.

Recall that ($1$-)equitability in  Section~\ref{sec:main} asserts that  given a matroid $\M=(E,\I)$ whose ground set $E$ can be partitioned into the union of $k$ disjoint bases, and given an arbitrary set $S\subseteq E$, there exists a partition $E=B_1\cup B_2\cup\ldots\cup B_k$ into $k$ bases such that $||S \cap B_i| - |S \cap B_j|| \leq 1$ for all $i,j \in [k]$. 
We now consider two disjoint sets $S_1,S_2\subseteq E$, and show that one can achieve the same guarantee for $S_1$, while having imbalance at most two for $S_2$; the bounds shown below are even slightly stronger.
In Example \ref{exm:tight}, we show that achieving $||S_i \cap B_j| - |S_i \cap B_{j'}|| \leq 1$ for all $j,j' \in [k]$ is not possible for both $i=1$ and $i=2$ at the same time. Hence, what we prove in Theorem \ref{thm:two-equitable} (and Corollary \ref{cor:two-equitable}), are the best one can hope for in the context of equally distributing two sets among disjoint bases. %

\begin{theorem}{\label{thm:two-equitable}}
Let $\M=(E,\I)$ be a matroid, and assume that the ground set $E$ can be partitioned into the union of $k$ disjoint bases. Let $S_1,S_2\subseteq E$ be two arbitrary disjoint sets. Then, there exists a partition $E=B_1\cup B_2\cup\ldots\cup B_k$ into $k$ bases such that the following hold:
\begin{enumerate}[label=(\roman*)]
\item\label{eq3}
$||S_1 \cap B_j| - |S_1 \cap B_{j'}|| \leq 1$ for all $j,j' \in [k]$; and %
\item\label{eq2}
$\left| |S_1 \cap B_j| - |S_1 \cap B_{j'}| \right|+\left| |S_2 \cap B_j| - |S_2 \cap B_{j'}| \right| \leq 2$ for all $j,j' \in [k]$; and
\item\label{eq1} 
  $||(S_1 \cup S_2) \cap B_j| - |(S_1 \cup S_2) \cap B_{j'}|| \leq \max_{i \in \{1,2\}}||S_i \cap B_j| - |S_i \cap B_{j'}||$ for all $j,j'\in [k]$.
\end{enumerate}
This implies, if $b^1_i = |B_i \cap S_1|$ and $b^2_i = |B_i \cap S_2|$, there exists $h,\ell \in \mathbb{Z}$ such that $\{(b^1_i,b^2_i) \mid i \in [n]\} \subseteq \{(h,\ell), (h, \ell+1), (h, \ell+2)\}$ or $\{(b^1_i,b^2_i) \mid i \in [n]\} \subseteq \{(h,\ell), (h,\ell+1), (h+1, \ell-1), (h+1,\ell)\}$. 
\end{theorem}
\begin{proof}
    For any partition $B=(B_1,B_2, \ldots, B_k)$ of $E$ into $k$ disjoint bases, let 
    \[
     \Phi(B)\coloneqq \sum_{j \in[k]} (|S_1 \cap B_j|^2+|S_2 \cap B_j|^2)\, ,\quad \quad    \Psi(B)\coloneqq \sum_{j \in [k]} |(S_1 \cup S_2) \cap B_j|^2\quad \mbox{and}\quad \Xi(B)\coloneqq\sum_{j \in[k]} |S_1 \cap B_j|^2.
    \]
    Now let $B=(B_1,B_2, \ldots, B_k)$ be a partition of $E$ into $k$ disjoint bases such that $\Phi(B)$ is minimum, subject to that, %
    $\Psi(B)$ is minimum, and subject to that, $\Xi(B)$ is minimum.
     We prove that $B$ satisfies conditions \ref{eq3}, \ref{eq2}, and \ref{eq1}.
  First we prove the following claim.
    \begin{claim}\label{claim1}
     For any $j,j' \in [k]$,
        if $|S_1 \cap B_j| > |S_1 \cap B_{j'}|$, then $|S_2 \cap B_j| \leq |S_2 \cap B_{j'}|$.
    \end{claim}
    \begin{proof}
    For simplicity of notation, let us assume that $j'=1$ and $j=2$. 
    Let $r=\mathrm{rk}(\M)$ denote the rank of the matroid. Thus, $|B_1|=|B_2|=r$ and $|E|=kr$.
    Let $a\coloneqq|S_1 \cap B_1|$, $b\coloneqq|S_1 \cap B_2|$, $c\coloneqq|S_2 \cap B_1|$, and $d\coloneqq|S_2 \cap B_2|$. Thus, we need to show that $b>a$ implies $d\le c$. Towards contradiction, assume $b>a$ and $d>c$. Let $S_3 \coloneqq (B_2 \cup B_1) \setminus (S_1 \cup S_2)$. We have 
    $|S_3 \cap B_2|<|S_3 \cap B_1|$ since $|S_3 \cap B_2| = r - b-d$ and $|S_3 \cap B_1| = r- a-c$. By Theorem \ref{thm:exchangable}, there exists $t \in B_2 \setminus S_3$ and $X \subseteq S_3$ such that $X$ is a $(t,S_3)$-exchangeable set. Let $\bar{B}_i \coloneqq B_i \Delta X$ for $i \in \{1,2\}$ and $\bar{B}_i\coloneqq B_i$ for all $i \in [k] \setminus \{1,2\}$. 
    Let us assume $t \in S_1$; the proof for $t\in S_2$ follows analogously.
    Then $|S_1 \cap \bar{B}_1| = a+1$, $|S_1 \cap \bar{B}_2|=b-1$, $|S_1 \cap \bar{B}_\ell| = |S_1 \cap B_\ell|$ for all $\ell \in [k] \setminus \{1,2\}$, and $|S_2 \cap \bar{B}_\ell| = |S_2 \cap B_\ell|$ for all $\ell \in [k]$. Thus, we get
    \[
    \begin{aligned}
        \Phi(\bar B) &= \Phi(B) - (a^2 + b^2) + ((a+1)^2 + (b-1)^2) \\
        &=\Phi(B) - 2(b-a-1) \\
        &\leq \Phi(B)\, ,\quad \mbox{and} \\
        \Psi(\bar B) &= \Psi(B) - (a+c)^2 - (b+d)^2 + (a+c+1)^2 + (b+d-1)^2 \\
        &=\Psi(B) - 2(b+d-a-c-1) \\
        &< \Psi(B)\, .
    \end{aligned}
    \]
    This contradicts the choice of $B$. Therefore, Claim \ref{claim1} holds. 
    \end{proof}

    Claim \ref{claim1} implies that condition \ref{eq1} holds.  We now turn to proving  \ref{eq3} and \ref{eq2}. Both conditions hold if $\max_{i \in \{1,2\},j,j' \in [k]}||S_i \cap B_j| - |S_i \cap B_{j'}|| \leq 1$. Let us assume that there exists $i \in \{1,2\}$ and $j,j' \in [k]$ such that $|S_i \cap B_j| - |S_i \cap B_{j'}| \geq 2$. Without loss of generality, let  $j=2$ and $j'=1$. Let $i'=3-i$. Let $a\coloneqq|S_i \cap B_1|$, $b\coloneqq |S_i \cap B_2|$, $c\coloneqq |S_{i'} \cap B_1|$, and $d\coloneqq |S_{i'} \cap B_2|$. 
    Thus, we are assuming $b-a\ge 2$. By Theorem \ref{thm:exchangable}, there exists $t \in B_1 \setminus S_i$ and $X \subseteq S_i$ such that $X$ is a $(t,S_i)$-exchangeable set. Let $\bar{B}_1 = B_1 \Delta X$, $\bar{B}_2 = B_2 \Delta X$, and $\bar{B}_j=B_j$ for all $j>2$. We have $|\bar{B}_1 \cap S_i| = a+1$ and $|\bar{B}_2 \cap S_i| = b-1$. By Claim \ref{claim1}, $c \geq d$. We have $b-a \geq 2$ and $c-d \geq 0$. We consider two cases.
    \paragraph{Case 1: $t \notin S_{i'}$.} We have $|S_{i'} \cap \bar{B}_1|=c$ and $|S_{i'} \cap \bar{B}_2|=d$. Therefore,
    \begin{align*}
        \Phi(\bar{B}) &= \Phi(B) - (a^2+b^2)+((a+1)^2 + (b-1)^2) \\
        &\leq \Phi(B) - 2(b-a-1) \\
        &< \Phi(B)\, ,
    \end{align*}
    which contradicts the choice of $B$. 
    \paragraph{Case 2: $t \in S_{\bar{i}}$.} We have $|S_{i'} \cap \bar{B}_1| =c-1$ and $|S_{i'} \cap \bar{B}_2|=d+1$.
    \begin{align*}
        \Phi(\bar{B}) &= \Phi(B) - (a^2+b^2+c^2+d^2)+((a+1)^2 + (b-1)^2 + (c-1)^2 + (d+1)^2) \\
        &\leq \Phi(B) - 2(b-a+c-d - 2)\, .
    \end{align*}
    If $b+c-a-d-2>0$, then $\Phi(\bar B) < \Phi(B)$ which is a contradiction. Otherwise, $b-a=2$, $c-d=0$, and hence $\Phi(B') = \Phi(B)$. If $i=2$, both conditions \ref{eq3} and \ref{eq2} hold. Thus, assume $i=1$ and $\bar{i}=2$. Since $t \in S_2$, we have $|(S_1 \cup S_2) \cap B_i| = |(S_1 \cup S_2) \cap B_{\ell}|$ for all $i, \ell \in [k]$. Thus, $\Psi(B') = \Psi(B)$. We have
    \begin{align*}
        \Xi(B') &= \sum_{j \in[k]} |S_1 \cap B'_j|^2 \\
        &= \sum_{j \in[k]} |S_1 \cap B_j|^2 - (a^2+b^2) + (a+1)^2 + (b-1)^2 \\
        &= \sum_{j \in[k]} |S_1 \cap B_j|^2 - 2(b-a-1) \\
        &< \sum_{j \in[k]} |S_1 \cap B_j|^2 \\
        &= \Xi(B),
    \end{align*}
    which again contradicts the choice of $B$.    
\end{proof}
Setting $k=2$, we get Corollary \ref{cor:two-equitable}. 
The next example shows that $2$-equitability is not true even for graphical matroids; hence, the above theorem gives in essence the best possible bound.

\begin{figure}
    \centering
    \includegraphics[width=0.15\linewidth]{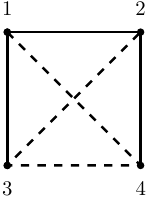}
    \caption{Solid edges correspond to the elements in $B_1$ and dashed edges to the elements in $B_2$. If $(C_1,C_2)$  equally divide $S_1 = \{(1,2), (3,4)\}$ and $S_2 = \{(1,3), (2,4)\}$, by symmetry we can assume $\{(1,2),(1,3)\} \subset C_1$ and $\{(2,4),(3,4)\} \subset C_2$. The set that includes $(2,3)$ has a cycle and therefore is not independent.}
    \label{fig:k4}
\end{figure}

\begin{example}\label{exm:tight}
    Consider $K_4$ with $V = \{1,2,3,4\}$. Let $S_1 = \{(1,2),(3,4)\}$ and $S_2 = \{(1,3),(2,4)\}$. See Figure \ref{fig:k4}. Note that $B_1 = \{(1,3),(1,2),(2,4)\}$ and $B_2 = \{(1,4),(2,3),(3,4)\}$ partition $E$ into two disjoint bases (i.e. spanning trees). However, there exists no partition of $E$ into two disjoint bases $B'_1,B'_2$ such that $|B'_1 \cap S_1| = 1$ and $|B'_1 \cap S_2|=1$. 
\end{example}

\section{Applications in Fair Division}\label{sec:fair-division}
In this section, we prove fairness guarantees in the matroid constrained setting when agents have bi-valued or tri-valued valuations. First in Lemma \ref{lem:shift-vals}, we prove that shifting the values by a constant preserves EF1 and MMS guarantees. 

\begin{lemma}\label{lem:shift-vals}
    Given a matroid constrained fair division instance $I=\langle N,E,(v_i)_{i\in N},\M\rangle$,  let $v'_i\, :\, 2^E \rightarrow \mathbb{R}_{\geq 0}$ be an additive valuation function such that for all $g \in E$, $v'_i(g)= v_i(g)-c_i$ for some constant $c_i$ for each $i\in N$. Then, 
    \begin{enumerate}[label=(\roman*)]
        \item if $A$ is an EF1 allocation for $I'=\langle N,E,(v'_i)_{i\in N},\M\rangle$, then $A$ is also an EF1 allocation for $I$, when $0\le c_i \leq \min_{g} v_i(g)$.
        \item if $A$ is an MMS allocation for $I'=\langle N,E,(v'_i)_{i\in N},\M\rangle$, then $A$ is also an MMS allocation for $I$.
    \end{enumerate}
    
\end{lemma}
\begin{proof}
    Let $r=\mathrm{rk}(\M)$ denote the rank of the matroid.  
    Let $A$ be an EF1 allocation for $I'$. Since $c_i \leq \min_{g} v_i(g)$, $I'$ consist of goods ($v_i(g) \geq 0$ for all $i \in N$ and $g \in E$) and the definition of EF1 is as in \Cref{def:ef1}.
    
    In order to prove $A$ is EF1 for $I$, it suffices to prove that for each $i\in A$, $v_i(A_i) \geq v_i(A_j - g)$ for some $g \in A_j$. There exists $g \in A_j$ such that 
    \begin{align*}
        v_i(A_i) &= v'_i(A_i) + c_i|A_i| \\
        &= v'_i(A_i) + c_i \cdot r \tag{$|A_i| = r$}\\
        &\geq v'_i(A_j - g) + c_i \cdot r\tag{$A$ is EF1 for $I'$} \\
        &= v_i(A_j - g) + c_i \tag{$|A_j - g| = r-1$}\\
        &\geq v_i(A_j - g). \tag{$c_i \geq 0$}
    \end{align*}
    
    Now, let $A$ be an MMS allocation for $I'$. In order to prove $A$ is MMS for $I$, it suffices to prove $v_i(A_i) \geq \mu^n_{v_i}(E)$. We have $\mu^n_{v_i}(E) = \mu^n_{v'_i}(E) + c_i \cdot r$. Also 
    \begin{align*}
        v_i(A_i) &= v'_i(A_i) + c_i \cdot r\\
        &\geq \mu^n_{v'_i}(E) + c_i \cdot r \\
        &= \mu^n_{v_i}(E).
    \end{align*}
\end{proof}

\trivalef*
\begin{proof}
    Let $v \coloneqq v_i$ for all $i \in N$. By Lemma \ref{lem:shift-vals}, we can assume that for all $g \in E$, $v(g) \in \{0,a,b\}$ for some $0<a<b$. Let 
    \[H \coloneqq \{g \mid v(g)=b\}\quad \mbox{and}\quad L \coloneqq \{g \mid v(g)=a\}\, .\]
     For any subset of items $T$, let $T^H \coloneq T \cap H$ and $T^L = T \cap L$ denote the ``high value'' and ``low value'' parts of $T$. Let $A$ be an allocation that equalizes $H$ and subject to that equalizes $L$ by setting $S_1=H$ and $S_2=L$ in Theorem \ref{thm:two-equitable}. Then $\{(|A_i^H|, |A_i^L|) \mid i \in N\} \subseteq \{(h,\ell), (h, \ell+1), (h, \ell+2)\}$ or $\{(|A_i^H|, |A_i^L|) \mid i \in N\} \subseteq \{(h,\ell), (h,\ell+1), (h+1, \ell-1), (h+1,\ell)\}$. In the second case, it is easy to check that the allocation is EF1. In the first case, if the allocation is not EF1 and $v_i(A_i) < v_i(A_j - g)$ for all $g \in A_j$, then $(|A_i^H|, |A_i^L|) = (h,\ell)$ and $(|A_j^H|, |A_j^L|) = (h,\ell+2)$. If $h=0$, by Theorem \ref{thm:main-equitable}, we can equalize $L$ up to one item. So assume $h>0$. Then, the EF1-envy implies
    \[
    a \cdot \ell + b \cdot h < a \cdot (\ell + 2) + b \cdot (h-1)\, .
\]
    Hence, $a > b/2$. Let $S \coloneqq (A_i \cup A_j) \setminus (H \cup L)$ denote the set of items in $A_i\cup A_j$ with $v(g)=0$.  
    Since $|S \cap A_j| < |S \cap A_i|$, by Theorem \ref{thm:exchangable}, there exists $t \in A_j \setminus S$ and $X \subseteq S+t$ such that $X$ is $(t,S)$-exchangeable. Let $\bar{A}_i \coloneqq A_i \Delta X$ and $\bar{A}_j \coloneqq A_j \Delta X$; we let $\bar{A_k}=A_k$ for all $k\in N\setminus\{i,j\}$. First we prove that there is no EF1-envy between agents $i$ and $j$ in $\bar{A}$.
    
    \paragraph{Case 1: \boldmath$t \in H$.} We have $(|\bar{A}_i^H|, |\bar{A}_i^L|) = (h+1,\ell)$ and $(|\bar{A}_j^H|, |\bar{A}_j^L|) = (h-1, \ell+2)$. Hence, we get 
    \begin{align*}
        v(\bar{A}_i) &= b \cdot (h+1) + a \cdot \ell \\
        &> b \cdot (h-1) + a \cdot (\ell+2) \tag{$b>a$} \\
        &= v(\bar{A}_j),
    \end{align*}
    and 
    \begin{align*}
        v(\bar{A}_j) &= b \cdot (h-1) + a \cdot (\ell+2) \\
        &> b \cdot h + a \cdot \ell \tag{$a > b/2$} \\
        &= v(\bar{A}_i - g). \tag{for all $g \in \bar{A}_i^H$}
    \end{align*}
    Therefore, there is no EF1-envy between agents $i$ and $j$ in $\bar{A}$.
    
    \paragraph{Case 2: \boldmath$t \in L$.} We have $(|\bar{A}_i^H|, |\bar{A}_i^L|) = (|\bar{A}_j^H|, |\bar{A}_j^L|) = (h,\ell+1)$. Clearly, there is no envy between agents $i$ and $j$ in $\bar{A}$.
       
    In the modified allocation, for all agents $k\in N$,
    \[ (|\bar{A}_k^H|, |\bar{A}_k^L|)\in \{(h-1,\ell+2), (h,\ell), (h, \ell+1), (h, \ell+2), (h+1, \ell)\}\, .
    \]
     It is easy to check that the only EF1-envy in this set can only be from an agent $p$ with $(|\bar{A}^H_p|,|\bar{A}^L_p|) = (h,\ell)$ to an agent $q$ with $(|\bar{A}^H_q|,|\bar{A}^L_q|) = (h,\ell+2)$. So as long as such a pair exists, repeat doing the above process. The final allocation is EF1.
\end{proof}

 Using Theorem \ref{thm:tri-val-ef1}, an EF1 allocation for two agents with (non-identical) tri-valued additive valuations can be obtained by a simple cut-and-choose protocol. Assuming that both of the agents have agent one's valuation function $v_1$, there exists a feasible EF1 allocation $(A_1,A_2)$. Now let agent two pick her favourite bundle among $A_1$ and $A_2$, and allocate the remaining bundle to agent one. The final allocation is clearly EF1.
\begin{corollary}[of Theorem \ref{thm:tri-val-ef1}]\label{cor:EF1-two-agents}
    Given a matroid constrained fair division instance $\langle N,E,(v_i)_{i\in N},\M\rangle$ with $|N|=2$, if the valuation $v_i$ is tri-valued additive for $i \in \{1,2\}$, there exists a feasible EF1 allocation.
\end{corollary}

\begin{lemma}\label{lem:reduction}
    Let $\langle N,E,(v_i)_{i\in N},\M\rangle$ be a matroid constrained fair division instance with binary additive valuations. Let $i \in N$ and $S=(S_1, \ldots, S_k, \ldots, S_n)$ be such that $S$ is a feasible allocation and $v_i(S_j) \leq \mu_{v_i}^n(E)$ for all $j \in [k]$. Then 
    \[\mu^{n-k}_{v_i}(E \setminus \bigcup_{j \in [k]}S_j) \geq \mu^n_{v_i}(E).\]
\end{lemma}
\begin{proof}
    We have
    \begin{align*}
        \mu^{n-k}_{v_i}\big(E \setminus \bigcup_{j \in [k]}S_j\big) &= \left\lfloor \frac{v_i(E \setminus \bigcup_{j \in [k]}S_j)}{n-k}\right\rfloor \tag{Corollary \ref{cor:mms-val}} \\   
        &= \left\lfloor \frac{v_i(E) - \sum_{j=1}^k v_i(S_j) }{n-k}\right\rfloor \\
        &\ge \left\lfloor\frac{v_i(E) - k \cdot \mu_i^n(E)}{n-k}\right\rfloor \tag{$v_i(S_j) \le \mu_i^n(E)$ for all $j \in [k]$}\\
        &= \left\lfloor\frac{v_i(E) - k \cdot \floor{v_i(E)/n}}{n-k}\right\rfloor \tag{Corollary \ref{cor:mms-val}} \\
        &\geq \left\lfloor\frac{v_i(E)} n\right\rfloor \\
        &= \mu^n_{v_i}(E). \tag{Corollary \ref{cor:mms-val}}
\end{align*}
\end{proof}

\mmsthm*
\begin{proof}
  By Lemma \ref{lem:shift-vals}, without loss of generality, we can assume for all agents $i$ and goods $g$, $v_i(g) \in \{0,a_i\}$. Since MMS is scale-invariant, without loss of generality we can assume the valuation functions are binary.    The rest of the proof is by induction on the number $n$ of agents. 
  
  For $n=1$, the theorem trivially holds. Assume for all $n'<n$ the claim holds. Fix an agent $a\in N$. Let $S=(S_1, \ldots, S_n)$ be an MMS partition of agent $a$. Create a bipartite graph $G$ with (nodes corresponding to) $S_1, \ldots, S_n$ on one side and (nodes corresponding to) agents $1,\ldots,n$ on the other side. There is an edge between $i$ and $S_k$ if and only if $v_i(S_k) \geq \mu_i^n(E)$. If $G$ has a perfect matching, then there exists an MMS allocation. Otherwise, without loss of generality let $T = \{S_1, \ldots, S_{k+1}\} \subseteq \{S_1, \ldots, S_n\}$ be a minimal Hall's violater set. It means that $|\Gamma(T)|<k+1$,  and for all $T' \subsetneq T$, $|\Gamma(T')| \geq |T'|$ where $\Gamma(K)$ is the neighborhood set of $K$ for all subsets $K$ of the nodes. Since $a$ is in the neighborhood of all non-empty subsets of $\{S_1, \ldots, S_n\}$, $k \geq 1$. By Hall's theorem, there exists a matching covering $S_1, \ldots, S_k$. Note that for all agents $i$ not covered by this matching, there is no edge between $i$ and $S_j$ for $j \in [k]$, otherwise, $|\Gamma(T)| \geq k+1$. By Lemma \ref{lem:reduction}, after removing the matching, i.e., assigning all agents in the matching the matched bundles, for all the remaining agents the MMS value cannot decrease. $(S_{k+1}, \ldots, S_n)$ is a partition of the remaining items into $n-k$ bases. Therefore, by the induction hypothesis, we can find an MMS allocation of the remaining goods to the remaining agents.
\end{proof}

\begin{remark}
    Note that in Theorem \ref{thm:MMS}, we do not require the valuations to be positive. So the theorem remains true even when agents have negative or mixed valuations. 
\end{remark}

\paragraph{Acknowledgements} We would like to thank Krist\'of B\'erczi and Tam\'as Schwarcz for their helpful comments and discussions. The second author would like to thank Siddharth Barman and Rohit Gurjar for introducing him to problems related to envy-freeness under matroid constraints.

\bibliographystyle{abbrv} 
\bibliography{Biblio}
\appendix
\section{Gabow’s Conjecture Implies the Existence of EF1
}
We now show that Gabow's conjecture (\Cref{conj:gabow}) implies the existence of matroid constrained EF1 allocations. 
Given a valuation function $v\, :\, E\to \R$, we let 
\[\Lex(A)\coloneqq \left(v(A_{\sigma(1)}),v(A_{\sigma(2)}),\dots, v(A_{\sigma(n)})\right)\, \]
 where $\sigma$  is a permutation such that $ v(A_{\sigma(1)})\leq v(A_{\sigma(2)})\leq \cdots \leq v(A_{\sigma(n)})$. %
\begin{lemma}\label{lem:GabowImpliesEF1}
Given a matroid constrained fair division instance $I=\langle N,E,(v_i)_{i\in N},\M\rangle$, let the valuations $v_i=v$ be identical, and additive. If \cref{conj:gabow} holds, then an allocation $A=(A_1,A_2,\ldots,A_n)$ of $n$ bases such that $\Lex(A)$ is lexicographically maximal is an EF1 allocation. 
\end{lemma}
\begin{proof}
    We first show this by contradiction when there are two agents. Let $A=(A_1=\{x_1,\dots,x_r\},A_2=\{y_1,\dots,y_r\})$ be one of the feasible allocations such that $\Lex(A)$ is lexicographically maximal. %
    Let $(x_1,\dots,x_r,$ $y_1,\dots,y_n)$ be the sequence guaranteed by \cref{conj:gabow} such that every set of $r$ cyclic consecutive elements forms a basis. Without loss of generality, let us assume that $v(A_1)<v(A_2).$ For $i\in [r]$, let $X_i=\{x_1,\dots , x_i\}$ and $Y_i=\{y_1,\dots , y_i\}$. Let $1\leq k \leq r$ be the smallest integer such that $v((A_1\setminus X_k)\cup Y_k)\geq v((A_2\setminus Y_k)\cup X_k)$. Since $\Lex(A)$ is lexicographically maximal,
    \[v((A_1\setminus X_{k-1})\cup Y_{k-1})\leq v(A_1) %
    \text{ and } v((A_2\setminus Y_k)\cup X_k)\leq v(A_1).\]
    Using these inequalities, 
    \begin{align*}
        v(E)-2v(A_1)&\leq v(E)-2v((A_2\setminus Y_k)\cup X_k) \\ 
        &= v(E)-2v((A_2\setminus Y_{k-1})\cup X_{k-1})+ 2(v(y_k)-v(x_k))\\
        &= v(E)-2\left(v(E)-v((A_1\setminus X_{k-1})\cup Y_{k-1}) \right) + 2(v(y_k)-v(x_k))\\
        &\leq 2v(A_1)-v(E) + 2(v(y_k)-v(x_k)), \end{align*}  implying
        $v(E)\le 2v(A_1)+v(y_k)-v(x_k)$.
    
    This further implies that 
\[v(A_2)-v(y_k)\leq v(A_1)-v(x_k)\, .\]
Since $v(x_k)\ge 0$ and $v(A_1)<v(A_2)$, this shows that the allocation $(A_1,A_2)$ is EF1.

    Now, we show this for arbitrary number of agents. Let  $A=(A_1,\dots,A_n)$ be an allocation with lexicographically maximal $\Lex(A)$. Without loss of generality, assume $v(A_1)\le v(A_2)\le \cdots \le v(A_n)$. Towards a contradiction, assume this allocation is not EF1. Thus, there must exist another bundle $A_k$ such that $(A_1,A_k)$ is not an EF1 allocation for agents $1$ and $k$. 
    Since $A$ maximizes $\Lex(A)$, the bases $A_1$ and $A_k$ must maximize the lexicographic valuation $\left(\min\{v(A_1),v(A_2)\},\max\{v(A_1),v(A_2)\}\right)$ on  $A_1\sqcup A_k$. By the above argument, this means that $A_1$ and $A_k$ form an EF1 allocation on the 2 agent instance on $A_1\sqcup A_k$, leading to a contradiction.
\end{proof}
\begin{remark}\em 
We note that the above proof requires Gabow's conjecture (\Cref{conj:gabow}). If we use the weaker statement of White's conjecture (\Cref{conj:white}) in the above argument, we  again obtain $v(A_2)-v(y)\leq v(A_1)-v(x)$, along with the assumption $v(A_2)>v(A_1)$. This implies EF1 as long as $y\in A_2$. However, in the case of a sequence as in White's conjecture, we do not necessarily have $y\in A_2$, and hence the current proof breaks (but the implication might still be true).
\end{remark}

\end{document}